\documentclass[11pt]{article}
\usepackage[english]{babel}
\usepackage{amssymb}
\usepackage{graphicx}
\usepackage{graphics}
\usepackage{amsthm}

\input pictex.tex

\newcommand{\clo}{\mathrm{S}^1}

\theoremstyle{definition}

\begin{document}

\date{}
\author{Andr\'es Navas}

\title{Group actions on 1-manifolds: 
a list of very concrete open questions}
\maketitle

\vspace{-0.1cm}

\noindent{\bf Abstract.} Over the last four decades, group actions on manifolds have deserved much attention by people coming from different fields, as for instance 
group theory, low-dimensional topology, foliation theory, functional analysis, and dynamical systems. This text focuses on actions on 1-manifolds. We present a (non exhaustive) 
list of very concrete open questions in the field, each of which is discussed in some detail and complemented with a large list of references, so that a clear panorama 
on the subject arises from the lecture.

\vspace{0.15cm}

\noindent{\bf Mathematics Subject Classification (2010):} 20F60, 22F50, 37B05, 37C85, 37E10, 57R30.

\vspace{0.15cm}

\noindent{\bf Keywords.} Group actions, dynamics, diffeomorphisms, 1-manifolds, circle.

\vspace{0.5cm}

From the very beginning, groups were recognized as mathematical objects endowed with a certain ``dynamics". For instance, Cayley realized every group 
as a group of permutations via left translations:
$$G \longrightarrow \mathcal{P} (G), \quad g \mapsto L_g: G \to G, \quad L_g (h) = gh.$$ 
For a finitely-generated group $G$, this action has a geometric realization: We can consider the so-called {\em Cayley graph} of $G$ whose vertices are the elements 
of $G$, two of which $f,g$ are relied by an edge whenever $g^{-1} f$ is a generator (or the inverse of a generator). Then $G$ becomes a subgroup of the group of 
automorphisms of this graph.

In general, the group of such automorphisms is larger than the group $G$. A classical result of Frucht \cite{frucht} consists on a slight modification of this construction 
so that the automorphisms group of the resulting graph actually coincides with $G$. In fact, there are uncountably many such modifications, even for the trivial group. 
The smallest nontrivial regular graph of degree 3 with trivial automorphism group is known as {\em the Frucht graph}, and is depicted below.

\begin{center}

\includegraphics[scale=0.4]{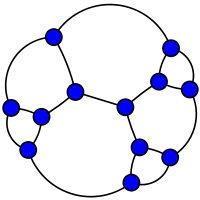}

Figure 1: The Frucht graph.

\end{center}

Another modification of Cayley's construction allows realizing every countable group as a group of homeomorphisms of the Cantor set. Assume for a while that such a 
group $G$ is infinite, and endow the space $M := \{ 0,1 \}^G = \{\varphi: G \to \{0,1\} \}$ with the product topology (metric). Then $M$ becomes a Cantor set, and $G$ 
faithfully acts on $M$ by shifting coordinates:
$L_g (\varphi) (h) := \varphi (g^{-1} h).$ 
Despite of its apparent simplicity, this shift action is fundamental in several contexts, 
and has attracted the attention of many people over the last decades \cite{boyle, cohen, H, popa2}.

In the case where $G$ is finite, one can modify the previous construction just by adding extra coordinates to the space $\{0,1\}^G$ on which the action 
is trivial. More interestingly, there is a single ``small" group of homeomorphisms of the Cantor set that contains all finite groups. To properly define it, we 
see the Cantor set as the boundary at infinite $\partial \Gamma$ of a regular tree $\Gamma$ of degree 3. Every proper, clopen ball in $\partial \Gamma$ 
can be canonically seen as the boundary at infinite of a rooted tree. We then consider the set of automorphisms of $\partial \Gamma$ that arise by cutting 
$\partial \Gamma$ into finitely many clopen balls and sending them into the pieces of another partition of $\partial \Gamma$ into clopen balls 
(with the same number of pieces), so that the restriction to each such ball is nothing but the canonical identification between these pieces 
viewed as boundaries of rooted trees. This yields the so-called {\em Thompson's group} $V$, which, among many remarkable properties, 
is finitely presented and simple. (See \cite{CFP} for more on this.) It is easy to see that $V$ contains all finite groups. 

Having realized every countable group as a group of homeomorphisms of a 0-dimensional space, one can ask whether some restriction arises when passing to higher 
dimension. Certainly, there are number of other motivations for considering this framework, perhaps the most transparent one coming from foliation theory. Indeed, to 
every group action of a finitely-generated group by homeomorphisms of a manifold $M$, one can associate a foliation by the classical procedure of suspension as 
follows: Letting $g_1,\ldots.g_k$ be a system of generators of $G$, we consider $S_k$, the surface of genus $k$, with fundamental group 
$$\pi_1 (S_k) = \left\langle a_1,\ldots, a_k, b_1, \ldots, b_k: \prod_{i=1}^{k} [a_i,b_i] = id \right\rangle.$$
In there, the generators $a_i$ are freely related, hence there is a homomorphism $\phi \!: \pi_1 (S_k) \to G$ sending $a_i$ into $g_i$ and $b_i$ into the identity. 
We then consider the product space $\Delta \times M$ endowed with the action of $\pi_1(S_k)$ given by $h (x,y) = (\bar{h}(x), \phi(h)(y))$, where $\bar{h}$ stands 
for the deck transformation on the Poincar\'e disc $\Delta$ associated to $h$. The quotient under this action is naturally a foliated, fibrated space with basis $S_k$ and 
fiber $M$, and the holonomy group of this foliation coincides with $G$. (See \cite{CC} for more on this construction.)

By the discussion above, and despite some remarkable recent progress \cite{b1,b2}, 
it seems impossible to develop a full theory of groups acting on manifolds. Here we restrict the discussion to the simplest case, namely, actions on 
1-dimensional spaces. In this context, the ordered structure of the phase space allows developing a very complete theory for actions by homeomorphisms, and the 
techniques coming from 1-dimensional dynamics allow the same for actions by diffeomorphisms. For each of such settings there are good references with very 
complete panoramas of the developments up to recent years: see \cite{order, ghys} and \cite{book}, respectively. This is the reason why we prefer to focus on 
challenging problems that remain unsolved, hoping that the reader will become motivated to work on some of them. 

\section{Actions of Kazhdan's groups}

In 1967, Kazhdan introduced a cohomological property and proved that it is satisfied by higher-rank simple Lie groups and their lattices, as for 
instance $\mathrm{SL}(n,\mathbb{Z})$ for $n \geq 3$ and their finite index subgroups \cite{kazhdan}. Since discrete groups satisfying this property 
are necessarily finitely generated, he proved finite generation for these lattices, thus solving a longstanding question. Since then, the so-called 
{\em Kazhdan's property (T)} has become one of the most important tools for studying actions and representations of Lie groups.

Although Kazhdan's original definition is somewhat technical, there is a more geometric property later introduced by Serre which turns 
out to be equivalent in the locally compact setting: a group satisfies {\em Serre's property (FH)} if every action by (affine) isometries on 
a Hilbert space has an invariant vector. (See \cite{BHV, margulis} for a full discussion on this.)

Property (T) has very strong consequences for the dynamics of group actions in different settings; see for instance \cite{bader,furman,popa,shalom, zimmer}. 
In what concerns actions on 1-dimensional spaces, a classical result pointing in this direction, due to Watatani and Alperin, states that every action of a 
group with property (T) by isometries of a real tree has a fixed point. In a more dynamical framework, Witte-Morris proved in \cite{witte} the 
following remarkable result: For $n \geq 3$, every action of a finite-index subgroup of $\mathrm{SL}(n, \mathbb{Z})$ by orientation-preserving 
homeomorphisms of the interval (resp. the circle) is trivial (resp. has finite image).

This theorem is even more remarkable because of its proof, which is amazingly elementary. However, it strongly relies on the existence 
of certain nilpotent subgroups inside the lattice, which do not arise in other ({\em e.g.} cocompact) cases. Despite some partial progress 
in this direction \cite{lucy1, lucy2,witte-2} (see \cite{witte-survey} for a full panorama on this), the following question remains open.

\vspace{0.25cm}

\noindent{\bf Question 1.} Does there exist a lattice in a higher-rank simple Lie group admitting a nontrivial action by orientation-preserving 
homeomorphisms of the interval ?

\vspace{0.25cm}

Notice that the statement above doesn't deal with actions on the circle. This is due to a theorem of Ghys \cite{ghys-invent}, which 
reduces the general case to that on the interval: Every action of a lattice in a higher-rank simple Lie group by orientation-preserving 
homeomorphisms of the circle has a finite orbit (hence a finite-index subgroup -which is still a lattice- fixes some interval). 

The question above can be rephrased in the more general setting of Kazhdan groups.

\vspace{0.25cm}

\noindent{\bf Question 2.} Does there exist an infinite, finitely-generated Kazhdan group of circle homeomorphisms ?

\vspace{0.25cm}

\noindent{\bf Question 3.} Does there exist a nontrivial (hence infinite) Kazhdan group of orientation-preserving homeomorphisms of the interval ?

\vspace{0.25cm}

A concrete result on this concerns actions by diffeomorphisms: If a finitely-generated group of $C^{3/2}$ circle diffeomorphisms satisfies property (T), 
then it is finite \cite{navas-kazhdan, book} (see \cite{yves} for the piecewise-smooth case). However, the situation is unclear in lower regularity. For instance, 
the group $G := \mathrm{SL} (2,\mathbb{Z}) \ltimes \mathbb{Z}^2$ has the {\em relative property (T)} (in the sense that for every action of $G$ by 
isometries of a Hilbert space, there is a vector that is invariant by $\mathbb{Z}^2$), yet it naturally embeds into the group of circle homeomorphisms. 
Indeed, the group $\mathrm{SL}(2,\mathbb{Z})$ acts projectively on the 2-fold covering of $\mathrm{S}^1$ -which is still a circle-, and blowing up an orbit 
one can easily insert an equivariant $\mathbb{Z}^2$-action. However, no action of this group is $C^1$ smoothable \cite{na-2, na-3}. 

The example above can be easily modified as follows: Letting $\mathbb{F}_2 \subset \mathrm{SL}(2,\mathbb{Z})$ be a finite-index subgroup, 
the semidirect product $G := \mathbb{F}_2 \ltimes \mathbb{Z}^2$ still has the relative property (T) (with respect to $\mathbb{Z}^2$). Moreover, 
starting with a free group of diffeomorphisms of the interval and using the blowing up procedure along a countable orbit, one can easily embed $G$ into 
the group of orientation-preserving homeomorphisms of the interval. Because of these examples, the answers to both Question 2 and 3 remain unclear.

There is a different, more dynamical approach to Question 1 above. Indeed, when dealing with the continuous case, actions on the interval and 
actions on the real line are equivalent. In the latter context, an easy argument shows that for every action of a finitely-generated group $G$ by 
orientation-preserving homeomorphisms of the real line without global fixed points, one of the following three possibilities occurs:

\vspace{0.1cm}

\noindent (i) There is a $\sigma$-finite measure $\mu$ that is invariant under the action.

\vspace{0.1cm}

\noindent (ii) The action is semiconjugate to a minimal action for which every small-enough interval is sent into a sequence of intervals that converge 
to a point under well-chosen group elements, but this property does not hold for every bounded interval. (Here, by a {\em semiconjugacy} we roughly 
mean a factor action for which the factor map is a continuous, non-decreasing, proper map of the real line.)

\vspace{0.1cm}

\noindent (iii) The action is semiconjugate to a minimal one for which the contraction property above holds for all bounded intervals.

\vspace{0.1cm}

Observe that a group may have actions of different type. (A good exercise is to build actions of $\mathbb{F}_2$ of each type.)  In case (i), the translation 
number homomorphism $g \mapsto \mu ([x,g(x)[)$ provides a nontrivial  homomorphism from $G$ into $\mathbb{R}$. 
In case (ii), it is not hard to see that, when looking at the minimal semiconjugate action, the map $\varphi$ that sends $x$ 
into the supremum of the points $y > x$ for which the interval $[x,y]$ can be contracted along group elements is an orientation-preserving homeomorphism 
of the real line that commutes with all elements of $G$ and satisfies $\varphi(x) > x$ for all $x$. Therefore, there is an induced $G$-action on the 
corresponding quotient space $\mathbb{R} / \!\sim$, where $x \sim \varphi(x)$, which is a topological circle. 

By the discussion above, case (i) cannot arise for infinite groups with property (T). As a direct consequence of Ghys' theorem stated above, 
case (ii) can neither arise for faithful actions of lattices in higher-rank simple Lie groups. Hence, if such a group admits an action (without 
global fixed points) on the real line, the action must satisfy property (iii).

\vspace{0.25cm}

\noindent{\bf Question 4.} Does there exist an infinite, finitely-generated group that acts on the real line all of whose actions 
by orientation-preserving homeomorphisms of the line without global fixed points are of type (iii) ?

\section{Cones and orders on groups}

Groups of orientation-preserving homeomorphisms of the real line are left orderable, that is, they admit total order relations that are invariant under left multiplication. 
Indeed, such a group can be ordered by prescribing a dense sequence $(x_n)$ of points in the line, and letting $f \prec g$ if the smallest $n$ for which 
$f(x_n) \!\neq\! g(x_n)$ is such that $f(x_n) \!<\! g(x_n)$. Conversely, for a countable left-orderable group, it is not hard to produce an action on the line. (See 
\cite{order, ghys} for more on this.) This may fail, however, for uncountable groups with cardinality equal to that of $\mathrm{Homeo}_+(\mathbb{R})$;  
see \cite{mann}.

The characterization above yields to a dynamical approach for the theory of left-orderable groups (which goes back to Dedekind and H\"older). In this 
view, a useful idea independently introduced by Ghys \cite{G} and Sikora \cite{sikora} consists in endowing the space $\mathcal{LO} (G)$ of all left-orders 
of a given left-orderable group $G$ with the Chabauty topology. (Two orders are close if they coincide over a ``large'' finite subset.) This provides a 
totally disconnected, compact space, which is metrizable in case $G$ is countable. (One can let $\mathrm{dist} (\prec,\prec') = 1/n$, where $n$ is 
the largest integer such that $\preceq$ and $\preceq'$ coincide over the set $A_n$ of a prescribed exhaustion of $G = \cup_i A_i$ by finite subsets.) 

A result of Linnell establishes that spaces of left orders are either finite or uncountable \cite{linnell-2} (see also \cite{CMR}). 
Left-orderable groups with finitely many orders were classified by Tararin \cite{order,km}: they are all solvable, the simplest examples being $\mathbb{Z}$ 
and the Klein bottle group $\langle a,b: bab = a \rangle$. In an opposite direction, for some classes of groups $G$, it is known that no left order is isolated 
in $\mathcal{LO}(G)$: solvable groups with infinitely many left orders \cite{RT}, free groups \cite{hs, ki, mc, na-4, rivas}, free products 
of groups \cite{rivas}, and surface groups \cite{ABR}. The following question remains, however, open.

\vspace{0.25cm}

\noindent{\bf Question 5.} Does there exist a finitely-generated, amenable, left-orderable group having an isolated order inside an infinite space of left orders~?

\vspace{0.25cm}

It is somewhat surprising that several classes of groups with infinitely many left orders do admit isolated left orders. Constructions have been proposed by  
different authors using quite distinct techniques: dynamical, group theoretical and combinatorial (see for instance \cite{dehornoy, dd, ito-1, ito-2, ito-3, matsu, na-is}). 
However, the most striking examples remain the first ones, namely, the braid groups $B_n$. To be more precise, let
$$B_n := \left\langle \sigma_1, \ldots, \sigma_{n-1} \!: \sigma_i \sigma_{i+1} \sigma_i = \sigma_{i+1} \sigma_i \sigma_{i+1}, \sigma_i\sigma_j=\sigma_j\sigma_i 
\mbox{ for } |i-j| \geq 2 \right\rangle$$
be the standard presentation of $B_n$. Denote $a_i := (\sigma_i \cdots \sigma_{n-1})^{(-1)^{i-1}}$, where $i \in \{1,\ldots,n-1\}$. Building on seminal work of 
Dehornoy \cite{deh-book}, Dubrovina and Dubrovin showed in \cite{dd} that $B_n$ admits the disjoint decomposition 
$$B_n = \left\langle a_1, \ldots, a_{n-1} \rangle^+ \cup \langle a_1^{-1}, \ldots, a_{n-1}^{-1} \right\rangle^+ \cup \{id\},$$
where $\langle \cdot \rangle^+$ stands for the semigroup generated by the corresponding set of elements. An easy argument then shows that the order 
whose elements larger than the identity are those in $\langle a_1, \ldots, a_{n-1} \rangle^+$ is well defined, total and left invariant; more importantly, it 
is isolated, since it is the only left order for which the elements $a_1,\ldots, a_{n-1}$ are all larger than the identity \cite{linnell-2}. 

Despite the apparent simplicity of the previous decomposition into finitely-generated positive and negative {\em cones}, 
no elementary proof is available. Finding an elementary approach is a challenging problem. The only nontrivial case 
that is well understood is that of $n=3$, where the decomposition is evident from the picture below. 

\begin{center}

\includegraphics[scale=0.42]{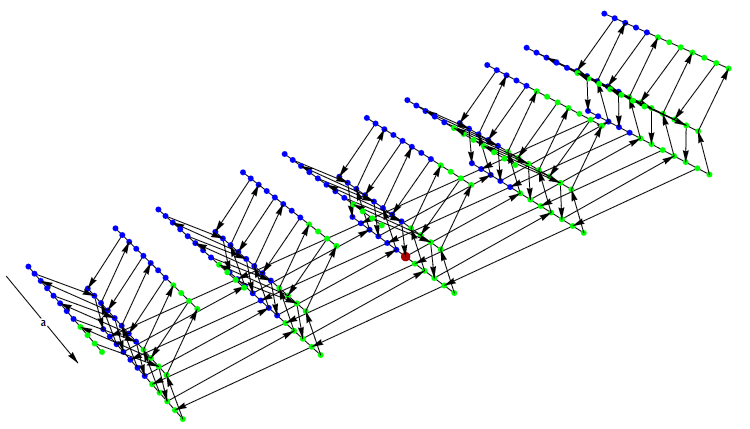}

Figure 2: The positive and negative cones of an isolated order on the Cayley graph of $B_3 = \langle a,b \!: ba^2b = a \rangle$ with respect to the generators 
$a := \sigma_1 \sigma_2$ and $b := \sigma_2^{-1}$.

\end{center}

Notice that a general left-orderable group $G$ acts on $\mathcal{LO} (G)$ by conjugacy: given a left order $\preceq$ and $ g \in G$, the conjugate of $\preceq$ under 
$g$ is the left order $\preceq_g$ for which $f_1 \preceq_g f_2$ if and only if $g f_1 g^{-1} \prec g f_2 g^{-1}$, which is equivalent to $f_1 g^{-1} \prec f_2 g^{-1}$.  

\vspace{0.25cm}

\noindent{\bf Question 6.} Does there exist a finitely-generated, left-orderable group for which the conjugacy action on its space of left orders is minimal 
(that is, all the orbits are dense) ?

\vspace{0.25cm}

It is not hard to show that free groups do admit left orders with a dense orbit under the conjugacy action \cite{mc, rivas}. However, this 
action is not minimal. Indeed,  free groups are bi-orderable (that is, they admit left orders that are also invariant under right multiplication), 
and a bi-order is, by definition, a fixed point for the conjugacy action.

The conjugacy action was brilliantly used by Witte-Morris in \cite{WW} to settle a question of Linnell \cite{linnell}, which was priorly raised -in the language of 
foliations- by Thurston \cite{thurston}: Every finitely-generated, left-orderable, amenable group admits a nontrivial homomorphism into $\mathbb{Z}$. Indeed,  
amenability provides an invariant probability measure for the conjugacy action. (Remind that one of the many definitions of {\em amenability} is that 
every action on a compact space admits an invariant probability measure.) Then the key idea is that, by the Poincar\'e recurrence theorem, 
left orders in the support of such a measure must satisfy a certain recurrence property, which reads as an algebraic property that is close 
to the Archimedean one (the so-called {\em Conradian property}; see \cite{NR-con}). This property allows obtaining the desired homomorphism.  
The next extension of the question (also proposed by Linnell), which is reminiscent of the Tits alternative, remains open.

\vspace{0.25cm}

\noindent{\bf Question 7.} Does there exist a finitely-generated, left-orderable group without free subgroups and admitting no nontrivial 
homomorphism into $\mathbb{Z}$~?

\vspace{0.25cm}

It is worth stressing that a negative answer to Question 4 above would imply a negative one to Question 7. Indeed, on the one hand, an action 
of type (i) provides a group homomorphism into $\mathbb{R}$ (via translation numbers), hence into $\mathbb{Z}$ for finitely-generated groups. 
On the other hand, as explained before, an action of type (ii) factors through a locally contracting action on the circle, which implies the presence
of a free subgroup by a theorem of Margulis \cite{margulis-cras} (see also \cite{ghys}).

A closely related question is the following.

\vspace{0.25cm}

\noindent{\bf Question 8.} Does there exist a finitely-generated, left-orderable group $G$ with no nontrivial homomorphism into $\mathbb{Z}$ 
and trivial group of bounded cohomology $\mathrm{H}^2_b (G, \mathbb{R})$~?

\vspace{0.25cm}

Again, a negative answer to Question 4 would also imply a negative one to this question. Indeed, locally contracting actions on $\mathrm{S}^1$ 
are parameterized (up to semiconjugacy) by a nontrivial cohomological class taking values in $\{ 0,1 \}$, according to a seminal work of Ghys \cite{ghys-mex} 
(see also \cite{ghys}). 

Besides these questions addressed for particular families of left-orderable groups, obstructions to left-orderability that go beyond torsion-freenes or 
the so-called {\em unique product property} UPP (namely, for each finite subset of the group, there is at least one product of two elements in this set that 
cannot be represented as another product of two elements in the set) are poorly understood. Although this goes beyond the scope of this text (and is one 
of the main lines of research of the theory of left-orderable groups), for specific families of dynamically defined groups, this has wide interest. A particular 
question on this was raised and nicely discussed by Calegari in \cite{cal-blog}.

\vspace{0.25cm}

\noindent{\bf Question 9.} Is the group of orientation-preserving homeomorphisms of the 2-disk that are the identity at the boundary left orderable ?

\vspace{0.25cm}

Notice that the group in discussion is torsion free, due to a classical result of Ker\'ekj\'art\'o \cite{ker} (see also \cite{kolev}). 
It is worth to stress that it is even unknown whether this group satisfies the UPP. 

\section{Groups of piecewise-projective homeomorphisms}

Groups of piecewise-affine homeomorphisms have been deeply studied in relation to Thompson's groups. Remind that Thompson's 
group $T$ is defined as the subgroup of (the previously introduced group) $V$ formed by the elements that respect the cyclic 
order in $\partial \Gamma$, the boundary of the homogeneous tree of degree-3. Besides, Thompson's group $F$ is the 
subgroup of $T$ formed by the elements that fix a specific point in $\partial \Gamma$ (say, the left-most point of the boundary 
of a clopen ball). Another realization arises in dimension 1: $T$ (resp. $F$) is a group of orientation-preserving, 
piecewise-affine homeomorphisms of the circle (rep. interval). Both groups are finitely presented \cite{CFP} and  have a dyadic 
nature, in the sense that the slopes of elements are integer powers of 2, and break points are 
dyadic rationals. One of the most challenging questions on these groups is the following.

\vspace{0.25cm}

\noindent{\bf Question 10.} Is Thompson's group $F$ amenable~?

\vspace{0.25cm}

Remind that a beautiful result of Brin and Squier establishes that the group of piecewise-affine homeomorphisms of the interval 
(hence $F$) doesn't contain free subgroups \cite{BS}. However, despite much effort by several people over the last decades (which 
includes several mistaken announcements pointing in the two possible directions), Question 10 remains as a kind of nightmare for the 
mathematical community; see \cite{ghys-agace}.

Besides the well-known question above, the algebraic structure of certain groups of piecewise-affine homeomorphisms, mainly 
generalizations of Thompson's groups \cite{stein}, is quite interesting. A concrete problem on them deals with distorted elements.
To properly state it, remind that, given a group $G$ with a finite generating system $\mathcal{G}$, the {\em word length} $\| g \|$ of $g \in G$ 
is the minimum number of factors needed to write $g$ as a product of elements in $\mathcal{G}$ (and their inverses). An element of 
infinite order $g \in G$ is said to be {\em distorted} if 
$$\lim_{n \to \infty} \frac{\| g^n \|}{n} = 0.$$
More generally, an element is said to be distorted in a general group whenever it is distorted inside a finitely-generated subgroup of this group.

Distorted elements naturally appear inside nilpotent groups, and have been extensively used to study rigidity phenomena of group actions 
on 2-manifolds \cite{CF,FH}. In the 1-dimensional setting, Avila proved that irrational rotations are distorted (in a very strong way) in the 
group of $C^{\infty}$ diffeomorphism \cite{avila}. However, despite some partial progress \cite{GL}, the following question is open.

\vspace{0.25cm}

\noindent{\bf Question 11.} Does the group of piecewise-affine circle homeomorphisms contain distorted elements~?

\vspace{0.25cm}

Beyond the piecewise-affine setting, the group of piecewise-projective homeomorphisms is a larger source of relevant examples 
of finitely-generated groups. In this direction, one could ask whether examples yielding to an affirmative answer to Question 3 may 
arise inside the group of piecewise-projective homeomorphisms of the line.

As a concrete example of an interesting group, remind that Thompson's group $T$ itself has a natural piecewise-projective, 
non piecewise-affine realization (which goes back 
to Thurston and, independently, to Ghys and Sergiescu): just replace dyadic rationals by rationals via the Minkowsky mark function, 
and change piecewise-affine maps by maps that are piecewise in $\mathrm{PSL}(2,\mathbb{Z})$. 

Among the new examples of groups constructed via this approach, the most remarkable is, with no doubt, the group $G_{LM}$ introduced by 
Lodha and Moore in \cite{LM}, which -as a group acting on the line- is generated by the homeomorphisms $f,g,h$ below (notice that $f,g$ 
generate a group isomorphic -actually, conjugate- to Thompson's group $F$):
$$f (t) : = t+1, 
\qquad 
g (t) := \left \{ \begin{array}{l}
t  \hspace{1.35cm} \mbox {if } t \leq 0,
\\ \\
\frac{t}{1-t} \hspace{0.95cm} \mbox{if } 0 \leq t \leq \frac12, 
\\ \\
3 - \frac{1}{t} \hspace{0.64cm} \mbox{if } \frac{1}{2} \leq t \leq 1,
\\ \\
t+1 \hspace{0.78cm} \mbox{if } t \geq 1,
\end{array} \right.
\quad \mbox{and } \quad
h ( t ) := \left \{ \begin{array}{l}
\frac{2t}{1+t} \hspace{1.2cm} \mbox {if } t \in [0,1],
\\ \\
t \hspace{1.64cm} \mbox{otherwise}. 
\end{array} \right.$$

Among its many remarkable properties, $G_{LM}$ has no free subgroup (by an easy extension of Brin-Squier's theorem mentioned above), 
it is non-amenable (due to prior work of Carri\`ere-Ghys \cite{CG} and Monod \cite{monod}), 
and has a finite presentation (this is the main technical contribution of \cite{LM}). 
Although it is not the first example of a group with these properties (see \cite{OS}), it has several 
other properties, as being torsion-free and of type F$_{\infty}$ (a property which is much stronger than being finitely presented; see \cite{lodha-1}). Recently, 
based on previous work relating smoothness with 1-dimensional hyperbolic dynamics in the solvable context \cite{BMNR}, this group was proven 
to be non $C^1$ smoothable in \cite{BLT} (see \cite{lodha} for a related result concerning a group closely related to $T$). This is to be compared 
with a classical result of Ghys and Sergiescu, according to which $T$ is topologically conjugate to a group of $C^{\infty}$ 
diffeomorphisms \cite{GS}. The following general question becomes natural.

\vspace{0.25cm}

\noindent{\bf Question 12.} What are the groups of piecewise-projective homeomorphisms of the interval/circle that are 
topologically conjugate to groups of $C^1$ diffeomorphisms~?

\section{The spectrum of sharp regularities for group actions}

Regularity issues appear as fundamental when dealing with both the dynamical properties of a given action and the algebraic constraints of 
the acting group. The original source of this goes back to Denjoy's classical theorem: Every $C^2$ circle diffeomorphism without periodic points 
is minimal. The $C^2$ hypothesis or, at least, a derivative with bounded variation, is crucial for this result. (The theorem is false in the $C^{1+\alpha}$ 
setting \cite{herman, tsuboi2}, and remains unknown for diffeomorphisms whose derivatives are $\tau$-continous with respect to the modulus of continuity 
$\tau (x) \!=\! | x \log(x)|$.) This is the reason why, when dealing with group actions (and, more generally, codimension-1 foliations), such an hypothesis is 
usually made. Nevertheless, in recent years, many new phenomena have been discovered in different regularities, thus enriching the theory.

One of the main problems to deal with in this direction is that of the optimal regularity. This problem is twofold. On the one hand, one looks for 
the maximal regularity that can be achieved, under topological conjugacy, of a given action. On the other hand, one asks for the maximal 
regularity in which a given group can faithfully act by varying the topological dynamics. A very concrete question in the latter 
direction is the following.

\vspace{0.25cm}

\noindent{\bf Question 13.} Given $0 < \alpha < \beta < 1$, does there exist a finitely-generated group of $C^{1+\alpha}$ diffeomorphisms 
of the circle/interval that does not embed into the group of $C^{1+\beta}$ diffeomorphisms~?

\vspace{0.25cm}

There are concrete reasons for restricting this problem only to regularities between $C^1$ and $C^2$. 
On the one hand, Kim and Koberda \cite{KK} have recently settled the analog of Question 13 for regularities larger than $C^2$, whereas 
the (discrete) Heisenberg group faithfully acts by $C^{1+\alpha}$ diffeomorphisms for any $\alpha < 1$ \cite{CJN}, but it does 
not embed into the group of $C^2$ diffeomorphisms \cite{PT}. On the other hand, Thurston 
gave the first examples of groups that are non $C^1$ smoothable via his remarkable stability theorem \cite{thurston} (see also \cite{BMNR, na-3}), 
while examples of groups of $C^1$ diffeomorphisms that are non $C^{1+\alpha}$ smoothable arise in relation to growth of groups \cite{gafa}. 
In an opposite direction, every countable group of circle homeomorphisms is topologically conjugate to a group of Lipschitz homeomorphisms, 
as it is shown below following the arguments of \cite{DKN-acta, compo}.

\vspace{0.25cm}

\noindent{\bf Example 1.} Let $G$ be a group with a finite, symmetric generating system $\mathcal{G}$ which acts by homeomorphisms 
of a compact 1-manifold $M$. Let $Leb$ denote the normalized Lebesgue measure on $M$. Given $\varepsilon > 0$, 
let $\bar{\mu}_{\epsilon}$ be the measure on $M$ defined as 
$$\bar{\mu}_{\varepsilon} := \sum_{f \in G} e^{-\varepsilon \| f \|} f_* (Leb),$$
where $\| f \|$ denotes the word-length of $g$ (with respect to $\mathcal{G}$). The measure $\bar{\mu}_{\varepsilon}$ has 
finite mass for $\varepsilon$ large enough. Indeed, 
$$\bar{\mu}_{\varepsilon} (M) = \sum_{n \geq 0} e^{-n \varepsilon} \big| S(n) \big|
\leq \sum_{n \geq 0} e^{-n \varepsilon} | \mathcal{G}| \big( |\mathcal{G}| - 1 \big)^{n-1}
= \frac{| \mathcal{G}|}{| \mathcal{G} | - 1} \sum_{n \geq 0} \left( \frac{|\mathcal{G}|-1}{e^{\varepsilon}} \right)^n,$$
where $| S(n) |$ stands for the cardinal of the set $S(n)$ of elements having word-length equal to $n$.
Moreover, since for every $g \in \mathcal{G}$ and all $f \in G$ it holds 
\hspace{0.01cm} $\| gf \| \leq \| f \|  + 1$, \hspace{0.01cm} we have 
\begin{equation}\label{quasi}
g_* (\bar{\mu}_{\varepsilon}) = \sum_{f \in G} e^{- \varepsilon \| f \| } (gf)_* (Leb) \leq
e^{\varepsilon} \sum_{f \in G} e^{-\varepsilon \| gf \|} (gf)_* (Leb)
= e^{\varepsilon} \bar{\mu}_{\varepsilon}.\end{equation}
Let $\mu_{\varepsilon}$ be the normalization of $\bar{\mu}_{\varepsilon}$. The probability measure $\mu_{\varepsilon}$ has total 
support and no atoms.  It is hence topologically equivalent to $Leb$ (dimension 1 is crucial here; see \cite{harrison,harrison2} for 
examples of non-smoothable homeomorphisms in higher dimension). By a change of coordinates sending 
$\mu_{\varepsilon}$ into $Leb$, relation (\ref{quasi}) becomes, for each interval $I \subset M$, 
$$|g^{-1}(I)| = g_* (Leb)(I) \leq e^{\varepsilon} Leb(I) = e^{\varepsilon} | I |.$$
This means that, in these new coordinates, $g^{-1}$ is Lipschitz with constant $\leq e^{\varepsilon}$.

\vspace{0.25cm}

It is interesting to specialize Question 12 to nilpotent group actions. Indeed, these actions are known to be $C^1$ smoothable \cite{FF,J,P}, 
though they are non $C^2$ smoothable unless the group is Abelian \cite{PT}. Moreover, the only settled case for Question 12 is that of a nilpotent 
group, namely, the group $G_4$ of $4 \times 4$ upper-triangular matrices with integer entries and 1's in the diagonal. In concrete terms, $G_4$ 
embeds into the group of $C^{1+\alpha}$ diffeomorphisms of the interval for every $\alpha < 1/2$, though it does not embed for 
$\alpha > 1/2$ \hspace{0.05cm} \cite{JNR}. (The case $\alpha \!=\! 1/2$ remains open; compare \cite{navas-israel}).

\section{Zero Lebesgue measure for exceptional minimal sets}

Differentiability issues are crucial in regard to ergodic type properties for actions. It follows from Denjoy's theorem quoted above that 
a single $C^2$ circle diffeomorphism cannot admit an {\em exceptional minimal set}, that is, a minimal invariant set homeomorphic to the Cantor set. 
Although the mathematical community took some time to realize that these sets may actually appear for group actions (the first explicit example appears 
in \cite{sa}), it is worth pointing out that these sets naturally arise, for instance, for Fuchsian groups (and, more generally, for groups with a Schottky 
dynamics), as well as for certain semiconjugates of Thompson's group $T$. The following question is due to Ghys and Sullivan.

\vspace{0.25cm}

\noindent{\bf Question 14.} Let $G$ be a finitely-generated group of $C^2$ circle diffeomorphisms. Assume that $G$ admits an exceptional minimal set 
$\Lambda$. Is the Lebesgue measure of $\Lambda$ equal to zero~?

\vspace{0.25cm}

An important recent progress towards the solution of this question is made in \cite{DKN-invent}, where it is answered in the affirmative for groups of 
real-analytic diffeomorphisms. Besides, an affirmative answer is provided, also in the real-analytic context, to another important question due to Hector.

\vspace{0.25cm}

\noindent{\bf Question 15.} Let $G$ be a finitely-generated group of $C^2$ circle diffeomorphisms. Assume that $G$ admits an exceptional minimal set 
$\Lambda$. Is the set of orbits of intervals of $ \mathrm{S}^1 \setminus \Lambda$ finite~?

\vspace{0.25cm}

Beyond having settled these two questions in the real-analytic context, the main contribution of \cite{DKN-invent} consists in proposing new ideas yielding 
to structure results for groups of circle diffeomorphisms admitting an exceptional minimal set (see \cite{DFKN} and the references therein for a full discussion 
on this). Indeed, so far, positive answers to these questions were known only in the expanding case (that is, whenever for every $x \!\in\! \Lambda$ there is 
$g \in G$ such that $Dg(x) > 1$; see \cite{duminy}), and for Markovian like dynamics \cite{CC1,CC2}. What is clear now is that, in the non expanding case, 
a certain Markovian structure must arise (see \cite{deroin} for a precise result in this direction in the conformal case). This view should also be useful 
to deal with the following classical question (conjecture) of Dippolito \cite{dippolito}.

\vspace{0.25cm}

\noindent{\bf Question 16.} Let $G$ be a finitely-generated group of $C^2$ circle diffeomorphisms. Assume that $G$ admits an exceptional minimal set 
$\Lambda$. Is the restriction of the action of $G$ to $\Lambda$ topologically conjugated to the action of a group of piecewise-affine homeomorphisms~? 

\vspace{0.25cm}

It should be pointed out that Questions 14, 15 and 16 have natural analogs for codimension-1 foliations. In this broader context, they all remain widely open, 
even in the (transversely) real-analytic setting. However, the ideas and techniques from \cite{DKN-invent} show that, also in this generality, structural issues 
are the right tools to deal with them.

\section{Ergodicity of minimal actions}

In case of minimal actions, a subtle issue concerns {\em ergodicity} (with respect to the Lebesgue measure), 
that is, the nonexistence of measurable invariant sets except for those having zero or full (Lebesgue) measure. 
The original motivation for this comes from a theorem independently proved by Katok \cite{katok} and Herman \cite{herman}: The 
action of a $C^2$ circle diffeomorphism without periodic points is ergodic (with respect to the Lebesgue measure). Notice that this result 
does not follow from Denjoy's  theorem (which only ensures minimality), since we know from the seminal work of Arnold \cite{arnold} 
that the (unique) invariant measure may be singular with respect to the Lebesgue measure.  

Katok's proof is performed via a classical {\em control of distortion} technique, which means that there is a uniform control on the ratio 
\hspace{0.015cm} $\sup Df_n /\inf Df_n$ \hspace{0.015cm} for the value of the derivatives on certain intervals along a well-chosen sequence 
of compositions $f_n$. This allows transferring geometric data from micro to macro scales, so that the proportion of the measures of different 
sets remains controlled when passing from one scale to another. Clearly, this avoids the existence of invariant sets of intermediate measure, 
thus proving ergodicity.

Herman-Katok's theorem deals with an ``elliptic'' context, whereas several classical ergodicity-like results (going back to Poincar\'e's linearization theorem) 
hold in an hyperbolic context. One hopes that a careful combination of both techniques would yield to an affirmative answer to the next question, 
also due to Ghys and Sullivan.

\vspace{0.25cm}

\noindent{\bf Question 17.} Let $G$ be a finitely-generated group of $C^2$ circle diffeomorphisms. If the action of $G$ is minimal, is it necessarily 
ergodic with respect to the Lebesgue measure~?

\vspace{0.25cm}

So far, an affirmative answer to this question is known in the case where the group is generated by elements that are $C^2$ close to rotations \cite{duminy}, 
for expanding actions \cite{duminy, DKN-moscow}, and for groups of real-analytic diffeomorphisms which are either free \cite{DKN-invent} or have infinitely 
many ends \cite{Jtop}. It is worth pointing out that the $C^2$ regularity hypothesis is crucial here; see for instance \cite{kodama}.

Again, Question 17 has a natural extension to the framework of codimension-1 foliations, where it remains widely open.

\section{Absolute continuity of the stationary measure}

Due to the absence of invariant measures for general groups actions, a useful tool to consider are the {\em stationary measures}, which correspond to 
probability measures that are invariant in mean. More precisely, given a probability distribution $p$ on a (say, finitely generated) group $G$ that acts by 
homeomorphisms of a compact metric space $M$, a probability measure $\mu$ on $M$ is said to be stationnary with respecto to $p$ if for every measurable 
subset $A \!\subset\! M$, one has
$$\mu (A) = \sum_{g \in G} p(g) \mu (g^{-1}(A)).$$ 
There are always stationary measures: this follows from a fixed-point argument or from Krylov-Bogoliuvob's argument consisting in taking means 
and passing to the limit in the (compact) space of probabilities on $M$. A crucial property is that, in case of uniqueness of the stationary measure (with 
respect to a given $p$), the action is ergodic with respect to this measure \cite{book}.

It is shown in \cite{antonov,DKN-acta} that, for a group of orientation-preserving circle homeomorphisms $G$ acting minimally, 
the stationary measure is unique with respect to each probability distribution $p$ on $G$ that is {\em non-degenerate} ({\em i.e.} 
the support of the measure generates $G$ as a semigroup). Hence, the next problem  becomes relevant in relation to Question 17 above.

\vspace{0.25cm}

\noindent{\bf Question 18.} Let $G$ be a finitely-generated group of $C^2$ orientation-preserving circle diffeomorphisms. Does there exist a 
non-degenerate probability distribution on $G$ for which the stationary measure is absolutely continuous with respect to the Lebesgue measure~?

\vspace{0.25cm}

A classical argument of ``balayage'' due to Furstenberg \cite{furst} solves this question for lattices in $\mathrm{PSL}(2,\mathbb{R})$. However, 
this strongly relies on the geometry of the Poincar\'e disk, and does not extend to general groups. 
Moreover, it is shown in \cite{DKN-moscow, guivarch} (see also \cite{peter}) 
that the resulting probability measure is singular with respect to the Lebesgue 
measure for non cocompact lattices whenever the distribution $p$ is {\em symmetric} ({\em i.e.} $p(g) = p(g^{-1})$ for all $g \in G$) and
finitely supported (and, more generally, for distributions with finite first moment). 
This also holds for groups with a Markovian dynamics for which there are {\em non-expandable points} ({\em i.e.} points $x$ such 
that $Dg(x) \leq 1$ for all $g \!\in\! G$), as for instance (the smooth realizations of) Thompson's group $T$. (Notice that, for the canonical action of 
$\mathrm{PSL} (2,\mathbb{Z})$, the point $[1 : 0]$ is non-expandable.)

Once again, Question 18 extends to the framework of codimension-1 foliations, where it remains widely open. (Uniqueness of the stationary measure in 
this setting and, more generally, in a transversely conformal framework, is the main content of \cite{DK-gafa}.)

A probability distribution on a group induces a random walk on it, many of whose properties reflect algebraic features of the group and translate into particular 
issues of the stationary measures. Remind that  to every probability distribution one can associate a ``maximal boundary'', which, roughly,  is a measurable 
space endowed with a ``contracting'' action having a unique stationary measure so that any other space of this type is a measurable factor of it \cite{furst2}. The study 
of this {\em Poisson-Furstenberg boundary} is one of the main topics in this area, and explicit computations are, in general, very hard \cite{erschler}. 
In our framework, a valuable result in this direction was obtained by Deroin, who proved in \cite{deroin-poisson} that for every group of 
smooth-enough circle diffeomorphisms with no finite orbit and whose action is locally discrete in a strong (and very precise)  sense, 
the Poisson-Furstenberg  boundary identifies with the circle endowed with the corresponding stationary measure provided the 
probability distribution on the group satisfies a certain finite-moment condition. Extending this result to more 
general groups is a challenging problem. In particular, the next question remains unsolved.

\vspace{0.25cm}

\noindent{\bf Question 19.} Given a symmetric, finitely-supported, non-degenerate probability distribution on Thompson's group $T$, does the 
Poisson-Furstenberg boundary of $T$ with respect to it identifies with the circle endowed with the corresponding stationary measure~?

\vspace{0.25cm}

Last but not least, random walks are also of interest for groups acting on the real line. In this setting, a nontrivial result is the existence of 
a (nonzero) $\sigma$-finite stationary measure for symmetric distributions on finitely-generated groups \cite{DKNP}. This is closely related 
to general recurrence type results for symmetric random walks on the line. One hopes that these ideas may be useful in dealing with 
Question 1, though no concrete result in this direction is known yet. 



 
\section{Structural stability and the space of representations}

Another aspect in which differentiability issues crucially appear concerns stability. Remind that, given positive numbers $r \leq s$, 
an action by $C^{s}$ diffeomorphisms is said to be $C^r$ structurally stable if every perturbation that is small enough in the $C^s$ 
topology is $C^r$ conjugate to it. (In the case $s = 0$, we allow semiconjugacies instead of conjugacies.) Usually, structural stability 
arises in hyperbolic contexts, and the situation in an elliptic type framework is less clear. The next question was formulated by 
Rosenberg more than 40 years ago (see for instance \cite{Ro}). 

\vspace{0.25cm}

\noindent{\bf Question 20.} Does there exist a faithful action of $\mathbb{Z}^2$ by $C^{\infty}$ orientation-preserving circle diffeomorphisms 
that is $C^{\infty}$ structurally stable ?

\vspace{0.25cm}

A closely related question, also due to Rosenberg, concerns the topology of the space of $\mathbb{Z}^2$-actions.

\vspace{0.25cm}

\noindent{\bf Question 21.} Given $r \geq 1$, is the subset of $\mathrm{Diff}_+^{r} (\mathrm{S}^1)^2$ consisting of pairs of commuting 
diffeomorphisms locally connected~?

\vspace{0.25cm}

These two questions have inspired very deep work of many people, including Herman and Yoccoz, 
who devoted their thesis to closely related problems.  However, despite all these efforts, they remain widely open. 
Among some recent progress concerning them, we can mention the proof of the connectedness of the space of commuting $C^{\infty}$ 
diffeomorphisms of the closed interval \cite{BE} (which, in its turn, has important consequences for codimension-1 foliations \cite{eynardcita}), 
and that of the path connectedness of the space of commuting $C^1$ diffeomorphisms 
of either the circle or the closed interval \cite{compo}. These two results apply in general to actions of $\mathbb{Z}^n$.

In a non-Abelian context, several other questions arise in relation to the structure of the space of actions. Among them, we can stress 
a single one concerning actions with an exceptional minimal set, for which the results from \cite{DKN-invent} point in a positive direction.

\vspace{0.25cm}

\noindent{\bf Question 22.} Given a faithful action $\phi_0$ of a finitely generated group $G$ by $C^{\infty}$ circle diffeomorphisms 
admitting an exceptional minimal set, does there exist a path $\phi_t$  of faithful actions of $G$ that is continuous in the $C^{\infty}$ 
topology and starts with $\phi_0 = \phi$  so that each $\phi_t$ admits an exceptional minimal set for $t < 1$ and $\phi_1$ is minimal~?

\vspace{0.25cm}

Quite surprisingly, structural stability is interesting even in the continuous setting. Indeed, the dictionary between left orders and actions on the interval 
shows that such an action is structurally stable if and only if a certain canonical left order arising from it is isolated in whole the space of left 
orders. Similarly, an action on the circle is structurally stable if and only if a natural ``cyclic order'' induced from it is an isolated point in the corresponding 
space of cyclic orders (endowed with the appropriate Chabauty topology; see \cite{MR}). In this regard, we may ask the following. (Compare Question 5.)

\vspace{0.25cm}

\noindent{\bf Question 23.} Let $G$ be a finitely-generated group of circle homeomorphisms whose action is $C^0$ structurally stable. 
Suppose that $G$ admits infinitely many non semiconjugate actions on the circle. Does $G$ contain a free subgroup in two generators~?

\section{Approximation by conjugacy and single diffeomorphisms}

Some of the connectedness results discussed above are obtained by constructing paths of conjugates of a given action. 
This idea is particularly simple and fruitful in very low regularity. We next give a quite elementary example to illustrate this.

\vspace{0.1cm}

\noindent{\bf Example 2.} As is well known, every circle homeomorphism has zero topological entropy. In most textbooks, this is proved by 
an easy counting argument of separated orbits for such an $f$. However, to show this, we can also 
follow the arguments of Example 1 for $G = \langle f \rangle$ 
and any $\varepsilon > 0$. Indeed, the outcome is that $f$ is topologically conjugate to a Lipschitz homeomorphism with Lipschitz constant 
$\leq e^{\varepsilon}$. By the invariance of entropy under topological 
conjugacy and its classical estimate in terms of the logarithm of the Lipschitz constant, we obtain that $h_{\mathrm{top}} (f) \leq \varepsilon$. 
Since this holds for all $\varepsilon > 0$, we must have $h_{\mathrm{top}} (f) = 0$.

\vspace{0.1cm}

The naive argument above still works for groups of subexponential growth, as for instance 
nilpotent groups \cite{compo}. Therefore, for all actions of these groups on 1-manifolds, the 
{\em geometric entropy} (as defined by Ghys, Langevin and Walczak in \cite{GLW}) always equals zero. More 
interestingly, a similar strategy should be useful to deal with groups of $C^1$ diffeomorphisms. Indeed, in this context, 
Hurder has shown in \cite{hurder} that zero entropy is a consequence of the absence of {\em resilient pairs}, which means 
that there are no elements $f,g$ such that \hspace{0.01cm} $x \!<\! f(x) \!<\! f(y) \!<\! g(x) \!<\! g(y) \!<\! y$ \hspace{0.01cm} for 
certain points $x,y$. (Notice that the converse holds even for homeomorphisms, as is follows from a classical counting argument.) The 
proof of  this fact is quite involved, and one hopes for an affirmative answer to the question below, which would immediately imply this result. 

\vspace{0.25cm}

\noindent{\bf Question 24.} Let $G$ be a finitely-generated group of $C^1$ circle diffeomorphisms. Suppose that 
$G$ has no resilient pairs. Given $\varepsilon > 0$, can $G$ be conjugated (by a homeomorphism) into a group 
of Lipschitz homeomorphisms for which the Lipschitz constants of the generators are all $\leq e^{\varepsilon}$~?
 
\vspace{0.25cm}
 
A particularly clarifying example on this concerns conjugates of $C^1$ diffeomorphisms without periodic points (that is, 
with irrational rotation number), as explained below. 

\vspace{0.25cm}

\noindent{\bf Example 3.} Remind that every cocycle $\varphi \!: M \to \mathbb{R}$ with respect to a continuous map $f \!: M \to M$ 
is cohomologous to each of its Birkhoff means. Indeed, letting
$$\psi_n := \frac{1}{n} \sum_{i=0}^{n-1} S_i \varphi, \qquad \mbox{ where } \qquad S_n \varphi:= \sum_{i=0}^{n-1} \varphi \circ f^i 
\quad \mbox{and} \quad S_0 \varphi := 0,$$
one easily checks the identity
$$\varphi - \frac{S_n \varphi}{n} = \psi_n - \psi_n \circ f.$$
If $f$ belongs to $\mathrm{Diff}_+^1 (\mathrm{S}^1)$, we can specialize this remark to $\varphi := \log (Df)$. 
Besides, if $f$ has irrational rotation number $\rho$, then an easy argument shows that $S_n (\log Df) / n \to 0$. Therefore,
$$\psi_n \circ f + \log Df - \psi_n \longrightarrow 0.$$
Adding a constant $c_n$ to $\psi_n$, we may assume that $\psi_n$ coincides with $\log D h_n$ for a $C^1$ diffeomorphism $h_n$. 
The relation above then becomes \hspace{0.1cm} 
$\log D (h_n f h_n^{-1}) \circ h_n \rightarrow 0,$ \hspace{0.1cm}
which shows that $h_n f h_n^{-1}$ converges to $R_{\rho}$ in the $C^1$ topology.

\vspace{0.25cm}

The argument above extends to actions of nilpotent groups, thus giving an affirmative answer to Question 25 for these groups \cite{compo}. 
However, this idea strongly uses the additive nature of the logarithm of the derivative, and it seems hard to directly extend it to higher regularity. 
Despite this, one expects that the use of a  Schwarzian-like derivative (cocycle) would yield to an affirmative answer to the following question.

\vspace{0.25cm}

\noindent{\bf Question 25.} Let $f$ be a $C^2$ circle diffeomorphism of irrational rotation number $\rho$. 
Does the set of $C^2$ conjugates of $f$ contain the rotation $R_{\rho}$ in its $C^2$-closure~?
 
\vspace{0.25cm}

The discussion above reveals that many natural questions still remain unsolved for single diffeomorphisms. Below we state two more of them.

\vspace{0.25cm}

\noindent{\bf Question 26.} For which values of $r > 1$ there exists $s \geq r$ such that for every $C^s$ circle diffeomorphism $f$ of irrational rotation number $\rho$, 
the sequence $f^{q_n}$ converges to the identity in the $C^r$ topology, where $(q_n)$ is the sequence of denominators of the rational approximations of $\rho$~?
 
\vspace{0.25cm}

This question is inspired by a fundamental result of Herman \cite{herman}, according to which one has the convergence \hspace{0.1cm } 
$f^{q_n} \!\to\! Id$ \hspace{0.1cm} in the $C^1$ topology  for $C^2$ circle diffeomorphisms of irrational rotation number (see also \cite{NT, yoccoz}). 
The answer to this question should consider a result of Yoccoz, who constructed a $C^{\infty}$ circle diffeomorphism with 
irrational rotation number and trivial centralizer \cite{yo-ast}. 

\vspace{0.25cm}

\noindent{\bf Question 27.} Let $f$ be a $C^2$ circle diffeomorphism of irrational rotation number $\rho$. Given $\varepsilon > 0$, let $M_{\varepsilon}$
 be the {\em mapping torus} of $f$ over $\mathrm{S}^1 \times [0,\varepsilon]$, that is, the surface obtained by identifying $(x,0) \sim (f(x), \varepsilon)$. 
 Let $\rho (\varepsilon) \in \mathbb{R} / \mathbb{Z}$ be such that $M_{\varepsilon}$ corresponds to the elliptic curve 
 $\mathbb{C} / ( \mathbb{Z} + i \rho(\varepsilon) \mathbb{Z})$. 
 Does $\rho(\varepsilon)$ converges to $\rho$ as $\varepsilon \to 0$~?
 
\vspace{0.25cm}

This question is due to Arnold. One hopes that recent progress on fine properties of circle diffeomorphisms should lead to a positive solution of it.

\section{Topological invariance of the Godbillon-Vey class}

The group of circle diffeomorphisms supports a remarkable cohomology class, 
namely, the {\em Godbillon-Vey class}, which is represented by the cocycle
$$(f,g) \mapsto \int_{\mathrm{S}^1} \log (Df) \hspace{0.03cm} D (\log D (g \circ f)).$$
Notice that, though this formula requires two derivatives, it can be naturally extended to $C^{3/2 + \varepsilon}$ diffeomorphisms (just pass half of the derivative 
from right to left; see \cite{HK, tsuboi-god}). However, no extension to $C^{1+\alpha}$ diffeomorphisms is possible for $\alpha$ small \cite{tsuboi1,tsuboi2}. 

According to a well-known result of Gelfand and Fuchs, the continuous cohomology of the whole group of $C^{\infty}$ circle diffeomorphisms 
is generated by two classes: the Euler class (which is the single generator in the $C^1$ setting), and the Godbillon-Vey class. Obviously, 
the Godbillon-Vey class induces (by restriction) a class in $H^2 (G,\mathbb{R})$ for every group $G$ of $C^{3/2+\varepsilon}$ circle 
diffeomorphisms. We refer to \cite{ghys-god, hurder-god} and the references therein for a panorama on this, including a full 
discussion on the next open question.

\vspace{0.25cm}

\noindent{\bf Question 28.} Is the (restriction of the) Godbillon-Vey class invariant under topological conjugacy for groups of $C^2$ diffeomorphisms~?

\vspace{0.25cm}

A first result in the positive direction was established by Raby, who proved invariance under conjugacy by $C^1$ diffeomorphisms \cite{raby}. Very 
soon after that, an alternative proof for this fact 
was proposed by Ghys and Tsuboi \cite{ghys-tsuboi}. Some years later, in \cite{HK}. Hurder and Katok proved invariance 
under conjugacies that are absolutely continuous (with an absolutely continuous inverse); see \cite{hilsum} for a recent result in the same direction. 

Ghys-Tsuboi's proof of Raby's theorem is of a dynamical nature. Indeed, in the most relevant cases of this framework, what it is proved is 
that $C^1$ conjugacies between groups of $C^r$ diffeomorphisms are automatically $C^r$ provided $r \geq 2$. This applies for instance 
to non-Abelian groups whose action is minimal. 

It is not hard to extend Ghys-Tsuboi's theorem to (bi-)Lipschitz conjugacies \cite{3-remarks, book}. However, absolutely continuous conjugacies 
are harder to deal with.

\vspace{0.25cm}

\noindent{\bf Question 29.} What are the groups of $C^2$ circle diffeomorphisms acting minimally for which the normalizer 
inside the group of absolutely continuous homeomorphisms coincides with that inside the group of diffeomorphisms~?

\section{On groups of real-analytic diffeomorphisms} 

The real-analytic framework offers new problems of wide interest even in the classical context. In this regard, remind that a celebrated result proved 
by Yoccoz in  \cite{yoccoz-cras} establishes that every real-analytic circle homeomorphism is minimal provided it has an irrational rotation number, thus 
extending Denjoy's theorem to this setting. (The same holds for $C^{\infty}$ homeomorphisms with non-flat singularities.) Ghys has asked 
whether this extends to the case where singularities may also arise for the inverse of the map.

\vspace{0.25cm}

\noindent{\bf Question 30.} Does Denjoy's theorem hold for circle homeomorphisms whose graphs are real-analytic~?

\vspace{0.25cm}

In a cohomological setting, another question concerns the validity of Geldfand-Fuch's theorem in the real-analytic case.

\vspace{0.25cm}

\noindent{\bf Question 31.}  Is the continuous cohomology of the group $\mathrm{Diff}^{\omega}_+(\clo)$ of orientation-preserving, real-analytic circle 
diffeomorphisms generated by the Euler and the Godbillon--Vey classes~?

\vspace{0.25cm}

A negative answer to this question would require the construction of a cocycle that uses real-analyticity in a crucial way. 
This would be somehow similar to Mather's homomorphism defined on the group of $C^1$ circle diffeomorphisms with 
derivatives having bounded variation. Remind that this is defined as 
$$f \mapsto \int_{\mathrm{S}^1} [ D (\log D f) ]_{\mathrm{reg}},$$
where $[ D (\log D f) ]_{\mathrm{reg}}$ stands for the regular part of the signed measure obtained as the derivative (in the sense of 
distributions) of the (finite total variation) function $\log D f$; see \cite{mather}. Such a homomorphism cannot exist in other regularities, 
because the corresponding groups of diffeomorphisms are known to be simple \cite{mather1,mather2}, except for class $C^2$ \cite{mather-raro}. 
By the way, though this is not related to real-analytic issues, this critical case must appear in any list of selected problems on the subject.

\vspace{0.25cm}

\noindent{\bf Question 32.} Is the group of orientation-preserving $C^2$ circle diffeomorphisms simple~?

\vspace{0.25cm}

Finally, we would like to focus on finitely-generated subgroups of $\mathrm{Diff}_+^{\omega} (\mathrm{S}^1)$. These have a tendency 
to exhibiting a much more rigid behavior than groups of diffeomorphisms. For example, though Thompson's groups act by $C^{\infty}$ 
diffeomorphisms, the group $F$ (hence $T$) does not faithfully act by real-analytic diffeomorphisms. One way to see this is by looking at 
solvable subgroups: $F$ contains such groups in arbitrary degree of solvability, though solvable groups of real-analytic diffeomorphisms 
of either the interval or  the circle are matabelian \cite{ghys-brasil}. (See however \cite{navas-solv} for algebraic constraints that apply to 
solvable groups of $C^2$ diffeomorphisms.)

Quite surprisingly, many algebraic issues that are known to hold or not to hold in the setting of $C^{\infty}$ diffeomorphisms are open in the 
real-analytic setting. For instance, it is unknown whether irrational rotations are distorted elements in $\mathrm{Diff}_+^{\omega} (\mathrm{S}^1)$. 
A more striking open question concerns the famous {\em Tits alternative}. 

\vspace{0.25cm}

\noindent{\bf Question 33.} Does the Tits alternative hold in $\mathrm{Diff}_+^{\omega} (\mathrm{S}^1)$~? More precisely, 
does every non-metabelian subgroup of this group contain a free subgroup~?

\vspace{0.25cm}

We refer to \cite{farb-shalen} for a partial result that reduces the general case to that of the interval. Notice that $F$ provides a negative answer to this 
question for groups of $C^{\infty}$ diffeomorphisms because of the aforementioned Ghys-Sergiescu's $C^{\infty}$ realization and Brin-Squier's theorem.

\vspace{0.3cm}

\noindent{\bf Added in Proof.} Question 13 has been recently solved by Kim and Koberda in a new version of \cite{KK}. 
I wish to thank Bassam Fayad and Sang-Hyun Kim for their remarks and corrections to an earlier version of this text.


\begin{footnotesize}


\vspace{0.25cm}

\noindent Andr\'es Navas (andres.navas@usach.cl)

\noindent Dpto de Matem\'atica y Ciencia de la Computaci\'on, Universidad de Santiago de Chile

\noindent Alameda 3363, Estaci\'on Central, Santiago, Chile

\end{footnotesize}

\end{document}